\documentclass[10pt]{article} 
\usepackage{amsmath, amssymb, amsthm, mathrsfs}
\usepackage{enumerate}
\usepackage{colordvi}
\usepackage{tikz}
\usepackage{graphicx}
\usepackage{hyperref}
\newtheorem{claim}{Claim}[section]
\newtheorem{theorem}[claim]{Theorem}

\newtheorem{remark}[claim]{Remark}
\newtheorem{corollary}[claim]{Corollary} 
\usepackage{color}

\definecolor{Myred}{cmyk}{0.0,1.0,1.0,0.00}
\definecolor{Mypurple}{rgb}{0.5,0.0,0.5}

\begin{document}

\title{Spectral estimates for Dirichlet Laplacian on spiral-shaped regions}

\author
{
Diana Barseghyan$^{1,2}$ and Pavel Exner$^{2,3}$
}
\date{
\small $^{1}$ Department of Mathematics, University of Ostrava,  30.dubna 22, 70103 Ostrava, Czech Republic \\
$^{2}$ Nuclear Physics Institute, Czech Academy of Sciences, 250 68 Rez, Czech Republic \\
$^{3}$Doppler Institute, Czech Technical University, B\v{r}ehov\'{a} 7, 11519 Prague, Czech Republic
\\
E-mails:\; diana.barseghyan@osu.cz, exner@ujf.cas.cz
}

\maketitle

\begin{abstract}
We derive spectral estimates of the Lieb-Thirring type for eigenvalues of Dirichlet Laplacians on strictly shrinking spiral-shaped domains.
\end{abstract}

\bigskip



\section{Introduction} \label{s: intro}

The dynamics of quantum particles confined to unbounded regions of various shapes is of interest not only from the physical point of view but also as a mathematical problem revealing interesting connections between spectral properties of the corresponding Hamiltonians and the confinement geometry. Problem of this type were discussed in numerous papers; for a survey and extensive bibliography we refer to the monograph \cite{EK15}. While some types of geometric perturbations such as bends or twists of straight tubes, both local and periodic, were investigated mathematically in considerable depth, some other mostly escaped attention. In particular, this is the case of spiral structures which appear in physics, for instance, as waveguides for cold atoms \cite{JLXZY15}. A mathematical analysis of Dirichlet Laplacians in spiral-shaped regions has been presented recently in \cite{ET21}; we refer to this paper also for references to other applications of spiral structures in different areas of physics.\footnote{Those include, in particular, electromagnetism and acoustics. Some of such systems are described by the Neumann Laplacian the spectral properties of which may be very different from the Dirichlet one as illustrated, e.g., in \cite{Si92}.}

Spiral regions are of many different type. A decisive factor for the spectral properties is the behavior of the coil width as we follow the spiral from the center to infinity; among those for which this quantity is monotonous, we can distinguish spirals expanding, asymptotically Archimedean, and shrinking -- definitions will be given below. In the present paper we are concerned with the latter type for which the Dirichlet Laplacian spectrum is purely discrete accumulating at the infinity. In such a situation it is natural to ask about moments of finite families of eigenvalues relative to a fixed energy value in the spirit of the estimates derived by Lieb and Thirring, Berezin, Lieb, and Li and Yau, cf.~\cite{LT76, Be72a, Be72b, Li73, LY83} and the monograph \cite{FLW22}. The cited results concern situations in which the motion is confined to finite region, or it is governed by a Schr\"odinger operator with a finite classically allowed region. The finiteness is not needed, however, similar estimates were also derived for infinite cusp-shaped regions \cite{GW11, BE13, EB14, BK19}.

Our main result, an estimate of the eigenvalue moments in terms of the geometric properties of the spiral, is presented and proved as Theorem~\ref{Main} in Sec.~\ref{Main result}. Before coming to it, we collect in the following section the necessary geometrical prerequisites. In particular, we will introduce locally orthogonal coordinates, sometimes called \emph{Fermi} or \emph{parallel}, that will allow us to rephrase the problem as spectral analysis of Dirichlet Laplacians on geometrically simpler cusped regions. In Sec.~\ref{s: conl}, we finish the paper with concluding remarks on the sharpness of the obtained bounds and on modifications of the result to the case of multi-armed spirals.

\section{Preliminaries} \label{s: prelim}

To begin with, let us describe the geometry of spiral-shaped regions. It is characterized be a curve $\Gamma$ which is the graph of an increasing function $r : \mathbb{R}_+ \to \mathbb{R}_+$ with $r(0)=0$, that is, the family of points $(r(\theta),\theta)$ in the polar coordinates. We note that spirals considered here are semi-infinite; a modification of our results to the case of `fully' infinite spirals, the example of which is the Simon's jellyroll mentioned above, is straightforward. The region we are interested in depends on the function $r$. Its closure is $\mathbb{R}^2$ provided that $\lim_{\theta\to\infty} r(\theta)= \infty$, in the opposite case it is the closed disc of radius $R:=\lim_{\theta\to\infty} r(\theta)$.

The assumed monotonicity of $r$ means that $\Gamma$ does not intersect itself which means, in particular, that the width function
\begin{equation}\label{a}
a(\theta):=\frac{1}{2\pi}\big(r(\theta)-r(\theta-2\pi)\big)
\end{equation}
is positive for any $\theta\ge 2\pi$. A spiral curve $\Gamma$ is called \emph{simple} if the corresponding $a(\cdot)$ is monotonous, and \emph{expanding} or \emph{shrinking} if this function is, respectively, increasing or decreasing in $\theta$ away from a neighborhood of the origin; these qualifications are labeled as \emph{strict} if $\lim_{\theta\to\infty}a(\theta)=\infty$ and $\lim_{\theta\to\infty}a(\theta)=0$ holds, respectively. A simple spiral laying between these two extremes, for which the limit is finite and nonzero, is called \emph{asymptotically Archimedean}.

The main object of our interest the two-dimensional Laplace operator with the Dirichlet condition imposed at the boundary represented by the curve $\Gamma$, in other words, the Dirichlet Laplacian $H_\Omega$ defined in the standard way \cite[Sec.~XIII.15]{RS78} on the open set $\Omega=\Omega_\Gamma=\mathbb{R}^2\setminus\Gamma$, or alternatively on $\mathcal{B}_R\setminus\Gamma$ if $\lim_{\theta\to\infty} r(\theta)<\infty$. As shown in \cite{ET21}, spectral properties of such operators on simple spiral regions depend strongly on the function $a(\cdot)$. For strictly expanding regions the spectrum is purely essential and covers the halfline $\mathbb{R}_+$. On the other hand, if the spiral $\Gamma$ is strictly shrinking the spectrum of $H_\Omega$ is purely discrete which is the situation we will be interested in.
\begin{remark}
{\rm In the intermediate case of asymptotically Archimedean spirals the spectrum may be more complicated. Its essential part covers the interval $\big[\frac14(\lim_{\theta\to\infty}a(\theta))^{-2}, \infty\big)$. The discrete part may be empty as in the case of the pure Archimedean spiral, but also infinite, accumulating at the threshold of $\sigma_\mathrm{ess}(H_\Omega)$ if the spiral is shrinking in the appropriate way \cite[Proposition~5.4]{ET21}. It is clear that neither of these situations allows us to derive bounds on eigenvalues moments analogous, say, to what one can derive in case of bent Dirichlet tubes \cite{ELW04}.
}
\end{remark}

A useful way to characterize the region $\Omega$, possibly with the exception of a neighborhood of the origin of the coordinates, is to employ the the Fermi (or parallel) coordinates, that is a locally orthogonal system in which the Cartesian coordinates of $\Gamma$ are written as
\begin{eqnarray}
\nonumber x_1(\theta, u) = r(\theta)\cos\theta-\frac{u}{\sqrt{\dot{r}(\theta)^2+r(\theta)^2}}
(\dot{r}(\theta) \sin\theta+r(\theta)\cos\theta),\\ \label{Fermi}
x_2(\theta, u)=r(\theta)\sin\theta+ \frac{u}{\sqrt{\dot{r}(\theta)^2+r(\theta)^2}}(\dot{r}(\theta)\cos\theta-r(\theta)\sin\theta),
\end{eqnarray}
where $u$ measures the distance of $(x_1, x_2)$ from $\Gamma$. A natural counterpart of the variable $u$ is the arc length of the spiral given by
\begin{equation}\label{s}
s(\theta)=\int_0^\theta\sqrt{\dot{r}(\theta)^2+r(\theta)^2}\,d\theta.
\end{equation}
We want to use relations \eqref{Fermi} and \eqref{s} to parametrize the region $\Omega$ with the coordinates $(s,u)$, possibly with the exception of a finite central part, as
\begin{equation}\label{parametrization}
\Omega_1=\Omega\cap\{(s, u): s>s_0\}= \big\{x(s, u): s >s_0,\: u \in(0, d(s))\big\}\,;
\end{equation}
here $s_0>0$ is a number depending on curve $\Gamma$ characterizing the excluded part, and $d(s)$ is the length of the inward normal starting from the point $x(s, 0)$ of $\Gamma$ towards the intersection with the previous coil of the spiral.

One more quantity associated with the spiral that we will need to state the result is its \emph{curvature} which is in terms of the angular variable given by
\begin{equation}\label{curvature}
\gamma(\theta)= \frac{r(\theta)^2+ 2\dot{r}(\theta)^2- r(\theta)\ddot{r}(\theta)}{
(r(\theta)^2+\dot{r}(\theta)^2)^{3/2}},
\end{equation}
provided, of course, that the derivatives make sense. Using the pull-back, $s\mapsto\theta(s)$, of the map \eqref{s} we can express it alternatively as a function of the arc length $s$, even if in general we lack an explicit expression; with abuse of notation we will write $\gamma(s)$ instead of $\gamma(\theta(s))$.

\section{Main result}\label{Main result}
\setcounter{equation}{0}

The domains to be considered are determined by the function $r$. In addition to its monotonicity and the requirement $r(0)=0$ we suppose that
 \begin{enumerate}[(a)]
 \setlength{\itemsep}{0pt}
\item $r$ is $C^2$-smooth and such that $\lim_{\theta\to\infty} \dot{r}(\theta)=0$; its second derivative is bounded and $\ddot{r}(\theta)<0$. \label{assa}
 \end{enumerate}
Under this assumption $\Gamma$ is shrinking, because in view of the concavity of $r$ the derivative $\dot{r}$ is decreasing and so is the function $\dot{a}(\theta):=\frac{1}{2\pi}\big(\dot{r}(\theta) -\dot{r}(\theta-2\pi)\big)$. Moreover, it is strictly shrinking, because $a(\theta) = \frac{1}{2\pi}\int_{\theta-2\pi}^{\theta} \dot{r}(\theta)\, \mathrm{d}\theta \to 0$ holds as $\theta\to\infty$, and furthermore, in view of \eqref{curvature} we have $\gamma(\theta)>0$. At the same time, the relation \eqref{curvature} in combination with the assumed boundedness of $\ddot{r}(\cdot)$ implies that $\gamma(\theta) = \mathcal{O}(r(\theta)^{-1})$ as $\theta\to\infty$, and since $d(s)\le 2\pi a(\theta(s))$ vanishes asymptotically, there is an $s_0>0$ such that $d(s)\gamma(s)<1$ holds for all $s\ge s_0$ which means, in particular, that in the corresponding part of $\Omega$ the Fermi coordinates are well defined.

The quantities we are interested in are moments of negative part of the operator $H_\Omega-\Lambda$ for a fixed energy $\Lambda$. We have the following bound:
\begin{theorem}\label{Main}
Let $\Omega=\Omega_\Gamma$ be a simple strictly shrinking domain determined by a spiral curve $\Gamma$ satisfying assumption \eqref{assa}, and let $H_\Omega$ be corresponding Dirichlet Laplacian. Then for any $\Lambda>0$ and $\sigma\ge\frac{3}{2}$ the following inequality holds,
\begin{equation}
\label{theorem}
\mathrm{tr}\left(H_\Omega-\Lambda\right)_-^\sigma \le\,\frac{L_{\sigma,1}
^{\mathrm{cl}}}{\pi} \big(\|W\|_\infty+\Lambda\big)^{\sigma+1}
\int_{\{d(s)\ge\pi(W(s)+\Lambda)^{-1/2}\}}d(s)\,\mathrm{d}s
+ c_1\Lambda^{\sigma+1}+c_2(\Lambda),
\end{equation}
where $c_1>0$ is an explicit constant given in \eqref{c} below, $c_2(\Lambda)$ is given in (\ref{cc}),
$\|\cdot\|_\infty:= \|\cdot\|_{L^\infty(s_0,\infty)}$, and
$L_{\sigma,1}^{\mathrm{cl}}$ is the semiclassical constant,
\begin{equation} \label{LTconstant}
L_{\sigma,1}^{\mathrm{cl}} := \frac{\Gamma(\sigma+1)}{\sqrt{4\pi} \Gamma(\sigma+\frac32)}\,;
 \end{equation}
the function $W$ in \eqref{theorem} is given by
\begin{equation}\label{W}
W(s):=\frac{\gamma^2(s)}{4(1-\gamma(s)d(s))^2}+\frac{d(s)|\ddot\gamma(s)|}{2(1-\gamma(s)d(s))^3}+\frac{5}{4}
\frac{d(s)^2|\dot\gamma(s)|^2}{(1-d(s)\gamma(s))^4}.
\end{equation}
Moreover, we have $$c_2(\Lambda)=\mathcal{O}
\left(\Lambda^2 \int_{\left\{d(s)\ge\frac{\pi}{\sqrt{\Lambda}}\right\}}d(s)\,ds\right)$$
for large values of $\Lambda$.
\end{theorem}
As a consequence, we get the asymptotic form of the bound:
\begin{corollary}\label{corollary}
In the regime $\Lambda\to\infty$, the inequality \eqref{theorem} becomes
\begin{equation}\label{asy_corollary}
\mathrm{tr}\left(H_\Omega-\Lambda\right)_-^\sigma \le\, \Lambda^{\sigma+1}\left(\frac{L_{\sigma,1}
^{\mathrm{cl}}}{\pi}\int_{\big\{d(s)\ge\frac{\pi}{\sqrt{\Lambda}}\big\}}d(s)\,\mathrm{d}s+c_1\right)\big(1+\bar{o}(1)\big).
\end{equation}
\end{corollary}

We note that in the asymptotic regime the curvature-related effects are suppressed being contained in the error term of \eqref{asy_corollary}. Beyond the asymptotics, they may be significant; for an analogy, recall the numerical result of \cite[Sec.~5.1]{ET21} about the number of eigenvalues in the Fermat spiral region which is at low energies larger than an estimate based on the spiral width only.

\medskip

\begin{proof}[Proof of Theorem (\ref{Main})]
As argued above, assumption \eqref{assa} makes it possible to use the Fermi coordinate parametrization (\ref{parametrization}) of $\Omega$ with some $s_0>0$. We employ the Neumann bracketing method \cite[Sec.~XIII.15]{RS78} which gives
\begin{equation}\label{bracketing}
H_\Omega\ge H_{\Omega_1}\oplus H_{\Omega_2},
\end{equation}
where $H_{\Omega_1}$ and $H_{\Omega_2}$ are the restrictions of $H_\Omega$ referring to the regions $\Omega_1\subset\Omega$ corresponding to the arc lengths $s>s_0$ and $\Omega_2:=\Omega\setminus\bar{\Omega}_1$, both having the additional Neumann condition at $s=s_0$.

Consider first the operator $H_{\Omega_1}$. According to \cite{ET21} the coordinates (\ref{Fermi}) allow us to pass from $H_{\Omega_1}$ to a unitarily equivalent operator $\widetilde{H}_{\widetilde{\Omega}_1}$ acting on the `straightened' region $\widetilde{\Omega}_1:=\{(s, u): s>s_0,\, 0<u<d(s)\}$ with Neumann boundary condition at $s=s_0$ and the Dirichlet condition on the rest of the boundary of $\Omega_1$ as follows,
\begin{equation}\label{unitary eq.}
\big(\widetilde{H}_{\widetilde{\Omega}_1}\psi\big)(s,u) = -\left(\frac{\partial}{\partial s}\frac{1}{(1-u\gamma(s))^2}\frac{\partial\psi}{\partial s}\right)(s,u)-\frac{\partial^2\psi}{\partial u^2}(s,u)+\widetilde{W}(s, u)\psi(s,u)\,,
\end{equation}
where
\begin{equation}\label{tildeW}
\widetilde{W}(s, u):=\frac{\gamma^2(s)}{4(1-u\gamma(s))^2}+\frac{u\ddot\gamma(s)}{2(1-u\gamma(s))^3}+\frac{5}{4}
\frac{u^2\dot\gamma(s)^2}{(1-u\gamma(s))^4}.
\end{equation}
It is straightforward to check that
\begin{equation}\label{est.}
\widetilde{H}_{\widetilde{\Omega}_1}\ge H_0,
\end{equation}
where operator $H_0=-\Delta-W$ acts on $L^2(\widetilde{\Omega}_1)$ and satisfies the same boundary conditions as $\widetilde{H}_{\widetilde{\Omega}_1}$, and $W$ is given by (\ref{W}).

Hence it is enough to deal with $H_0$.

We take inspiration from \cite{LW00, W08} and use a variational argument to reduce the problem to a Lieb-Thirring inequality with an operator-valued potential. Given a function $g\in\,C^\infty (\widetilde{\Omega}_1)$ with zero trace at the `transverse part' of the boundary, that is, at the points $\{s\in(s_0, \infty),\, u=0,\,d(s)\}$, and a number $\Lambda>0$, we can express the value of the corresponding quadratic form as

\begin{align*}
& \int_{\widetilde{\Omega}_1}\left(|\nabla\,g(s,u)|^2-(W+\Lambda)(s)|g(s,u)|^2\right)\,\mathrm{d}s\,\mathrm{d}u \\ & =\int_{\widetilde{\Omega}_1}\Big|\frac{\partial g}{\partial s}(s,u)\Big|^2\,\mathrm{d}s\,\mathrm{d}u
+\int_{s_0}^\infty\,\mathrm{d}s \int_0^{d(s)}\bigg(\Big|\frac{\partial g}{\partial u}(s,u)\Big|^2
- (W(s)+\Lambda) |g(s,u)|^2\bigg)\,\mathrm{d}u \\
& =\int_{\widetilde{\Omega}_1}\Big|\frac{\partial g}{\partial s}(s,u)\Big|^2\,\mathrm{d}s\,\mathrm{d}u +\int_{s_0}^\infty\left\langle\,L(s,W,\Lambda)g(s,\cdot), g(s,\cdot)\right\rangle_{L^2(0, d(s))}
\,\mathrm{d}s\,,
\end{align*}
where $L(s,W,\Lambda)$ is the Sturm-Liouville operator
$$
L(s, W, \Lambda)=-\frac{\mathrm{d}^2}{\mathrm{d}u^2} -W(s)-\Lambda
$$
defined on $L^2(0, d(s))$ with Dirichlet conditions at $u=0$ and $u=d(s)$.

Next we consider the complement of $\widetilde{\Omega}_1$ to the halfplane $\{s\ge s_0, u\in\mathbb{R}\}$ and denote its interior as $\widetilde{\Omega}_\mathrm{c}$. We take arbitrary functions $g\in C^\infty(\widetilde{\Omega}_1)$ and $v\in C^\infty(\widetilde{\Omega}_\mathrm{c})$, both having zero trace at $\{s\in(s_0, \infty),\, u=0,\,d(s)\}$; extending them by zero to the complements of $\widetilde{\Omega}_1$ and $\widetilde{\Omega}_\mathrm{c}$, respectively, we can regard them as functions in the whole halfplane. Similarly we extend $L(s,W,\Lambda)$ to the operator on $L^2(\mathbb{R})$ acting as $L(s,W,\Lambda)\oplus 0$ with the zero component on $\mathbb{R}\setminus[0, d(s)]$. For their sum, $h=g+v$, we then have
\begin{align*}
& \|\nabla\,g\|^2_{L^2(\widetilde{\Omega}_1)} +\|\nabla\,v\|^2_{L^2(\widetilde{\Omega}_2)} -\int_{\widetilde{\Omega}_1}(W(s)+\Lambda) |g(s,u)|^2\,\mathrm{d}s\,\mathrm{d}u
\\ & \quad
\ge\int_{\{s>s_0, u\in\mathbb{R}\}}\Big|\frac{\partial h}{\partial s}(s,u)\Big|^2\,\mathrm{d}s\,\mathrm{d} u +\int_{s_0}^\infty\left\langle
\,L(s,W,\Lambda)\,h(s,\cdot),\,h(s,\cdot)\right\rangle_{L^2(\mathbb{R})}
\,\mathrm{d}s\,.
\end{align*}
The left-hand side of this inequality is the quadratic form corresponding to the direct sum of operator $H_0-\Lambda$ and the Laplace operator defined on $\widetilde{\Omega}_\mathrm{c}$ with the Neumann boundary conditions at $s=s_0$ and Dirichlet conditions at the rest part of the boundary, while the right-hand side is the form associated with the operator
$$
-\frac{\partial^2}{\partial\,s^2}\otimes\,I_{L^2(\mathbb{R})} +L(s,W,\Lambda)\,,
$$
the form domain of which is larger, namely $\mathcal{H}^1\left((s_0, \infty), L^2(\mathbb{R})\right)$. Since the Laplace operator is positive, the minimax principle allows us to infer that
\begin{equation}\label{op.val.}
\mathrm{tr}\,(H_0-\Lambda)_-^\sigma \le \,\mathrm{tr}\Big(-\frac{\partial^2}{\partial\,s^2}\otimes \,I_{L^2(\mathbb{R})}+ L(s,W,\Lambda)\Big)_-^\sigma
\end{equation}
holds for any nonnegative number $\sigma$. This makes it possible to employ the following version of Lieb-Thirring inequality for operator-valued potentials (the proof of which is given in Appendix):
\begin{theorem}\label{operator valued1}
Let $Q(s),\,s\ge0$, be a family of self-adjoint operators on $L^2(\mathbb{R})$ with discrete spectrum. Then the following estimate holds for any $\sigma\ge 3/2$,
\begin{equation}\label{new result}
\mathrm{tr}\Big(-\frac{\partial^2}{\partial\,s^2}\otimes \,I_{L^2(\mathbb{R})}+ Q(s)\Big)_-^\sigma
\le\, L_{\sigma, 1}^{\mathrm{cl}}\int_0^\infty\mathrm{tr}\,Q_-^{\sigma+1/2}(s)\,\mathrm{d}s+\frac12 N\lambda_1^{(N)},
\end{equation}
where the operator on the left-hand side acts on the space $L^2(\mathbb{R}_+, L^2(\mathbb{R}))$ of vector-valued functions with the Neumann boundary condition at $s=0$, the number of its negative eigenvalues is denoted by $N$, and
\begin{equation}\label{Nlambda}
-\lambda_1^{(N)}:=\mathrm{inf}\Big(\,\sigma\Big(-\frac{\partial^2}{\partial\,s^2}\otimes\,I+P_N Q(s) P_N\Big)\cap (-\infty, 0)\Big),
\end{equation}
where $P_N$ is the projection on the span of the eigenfunctions referring to the negative eigenvalues of the operator $-\frac{\partial^2}{\partial\,s^2}\otimes \,I_{L^2(\mathbb{R})}+ Q(s)$.
\end{theorem}
\begin{remark}\label{s0}
{\rm One can easily check that Theorem(\ref{operator valued1}) holds true if $L^2(\mathbb{R}_+, L^2(\mathbb{R}))$ is replaced by $L^2((s_0, \infty), L^2(\mathbb{R}))$ with some  $s_0\in\mathbb{R}$; the only price to pay is the change of the integration interval on the right-hand side of \eqref{new result} to $(s_0, \infty)$.}
\end{remark}

Let us now return to the estimate \eqref{op.val.}. In order to apply Theorem~\ref{operator valued1} we have to find a suitable upper bound for the number of negative eigenvalues of $-\frac{\partial^2}{\partial\,s^2}\otimes \,I_{L^2(\mathbb{R})}+ L(s,W,\Lambda)$. We employ a `mirroring' trick the idea of which belongs to Rupert Frank \cite{F22}. It consists of introducing extended potential functions,
\begin{gather}\label{notation1}
\hat{W}(s_1+t)=W(s_1+t) \quad \text{and}\quad \hat{d}(s_1+t)=d(s_1+t)\quad \text{if}\quad t \ge0,\\[.2em]
\label{notation2} \hat{W}(s_1+t)=W(s_1-t)\quad \text{and} \quad \hat{d}(s_1+t)=d(s_1-t) \quad \text{if}\quad t<0,
\end{gather}
with which it is easy to see that for any $\Lambda'>0$ the spectrum of operator $\hat{H}=-\frac{\partial^2}{\partial\,s^2}\otimes \,I_{L^2(\mathbb{R})}+ L(s,\hat{W},\Lambda')$ defined on $L^2(\mathbb{R}, L^2(\mathbb{R}))$ contains the spectrum of operator of it `Neumann half', $-\frac{\partial^2}{\partial\,s^2}\otimes \,I_{L^2(\mathbb{R})}+ L(s,W,\Lambda')$. This implies
$$
\mathrm{tr}\Big(-\frac{\partial^2}{\partial\,s^2}\otimes \,I_{L^2(\mathbb{R})}+ L(s,W,\Lambda')\Big)_-^\sigma\le \mathrm{tr}\big(\hat{H}\big)_-^\sigma,\quad \sigma\ge0.
$$
Next we use to fact the validity of spectral estimates of this type can be extended to smaller values of the power $\sigma$ at the price of having a multiplicative coefficient $r(\sigma,1)$ on the right hand side, cf.~inequality \eqref{ELW} in Remark~\ref{remark} below in which also the explicit knowledge of eigenvalues of the transverse part of the operator is used. From that inequality in combination with above estimate  we get for $\sigma=1/2$ and $\Lambda'=2\Lambda$ the bound
\begin{align*}
& \mathrm{tr}\,\Big(-\frac{\partial^2}{\partial\,s^2}\otimes \,I_{L^2(\mathbb{R})}+ L(s,W,2\Lambda)\Big)_-^{1/2}\\[.3em]
& \le\,\frac{2r(1/2,1)\,L_{1/2,1} ^{\mathrm{cl}}}{\pi} \big(\|W\|_{L^\infty(s_0, \infty)}+2\Lambda\big)^{3/2} \int_{\big\{d(s)\ge\pi(W(s)+2\Lambda)^{-1/2}\big\}} d(s)\,\mathrm{d}s
\end{align*}
Let $N(\Lambda)$ be the number of negative eigenvalues of $-\frac{\partial^2}{\partial\,s^2}\otimes \,I_{L^2(\mathbb{R})}+ L(s,W,\Lambda)$. Using the last inequality together with a simple estimate,
\begin{align*}
& \mathrm{tr}\,\Big(-\frac{\partial^2}{\partial\,s^2}\otimes \,I_{L^2(\mathbb{R})}+ L(s,W,2\Lambda)\Big)_-^{1/2} \\[.3em] & =\mathrm{tr}\,\Big(-\frac{\partial^2}{\partial\,s^2}\otimes \,I_{L^2(\mathbb{R})}+ L(s,W,\Lambda)-\Lambda\Big)_-^{1/2} \ge \sqrt{\Lambda}\, N(\Lambda),
\end{align*}
we infer that
\begin{equation}\label{count.}
N(\Lambda)\le\, \frac{2r(1/2,1)\,L_{1/2,1} ^{\mathrm{cl}}}{\pi} \big(\|W\|_{L^\infty(s_0, \infty)}+2\Lambda\big)^{3/2} \int_{\big\{d(s)\ge\pi(W(s)+2\Lambda)^{-1/2}\big\}} d(s)\,\mathrm{d}s
\end{equation}

Now we are in position to apply the operator-valued version of Lieb-Thirring inequality. Let $P_{N(\Lambda)}$ be the projection on the span of the eigenfunctions corresponding to the negative eigenvalues of operator $-\frac{\partial^2}{\partial\,s^2}\otimes \,I_{L^2(\mathbb{R})}+ L(s,W,\Lambda)$. Combining Theorem~\ref{operator valued1} with \eqref{count.} we arrive at the estimate
\begin{equation}\label{op.val.new}
\mathrm{tr}\,\Big(-\frac{\partial^2}{\partial\,s^2}\otimes \,I_{L^2(\mathbb{R})}+ L(s,W,\Lambda)\Big)_-^\sigma \le \,L_{\sigma,1}^{\mathrm{cl}} \int_0^\infty\mathrm{tr}\,L(s,W,\Lambda)_-^{\sigma+1/2}\, \mathrm{d}s+c_2(\Lambda),
\end{equation}
where
\begin{equation}\label{cc}
c_2(\Lambda):=\frac{r(1/2,1)\,L_{1/2,1}^{\mathrm{cl}}\lambda_1(\Lambda)}{\pi \sqrt{\Lambda}} \big(\|W\|_{L^\infty(s_0, \infty)}+2\Lambda\big)^{3/2} \int_{\big\{d(s)\ge\pi(W(s)+2\Lambda)^{-1/2}\big\}} d(s)\,\mathrm{d}s
\end{equation}
and
\begin{equation}\label{ccc}
-\lambda_1(\Lambda):=\mathrm{inf}\,\Big(\sigma\Big(-\frac{\partial^2}{\partial\,s^2}\otimes \,I_{L^2(\mathbb{R})}+P_{N(\Lambda)} L(s,W,\Lambda) P_{N(\Lambda)}\Big)\cap(-\infty, 0)\Big).
\end{equation}

To make use of \eqref{op.val.new} we need to know the negative eigenvalues of $L(s,W,\Lambda)$. Those, however, are easy to be found using the fact $L(s,0,0)$ is the Dirichlet Laplacian on $(0, d(s))$ with the eigenvalues $\big(\frac{\pi j}{d(s)}\big)^2,\: j=1,2,\dots\,,$ and the potential is independent of the transverse variable $u$. Consequently, the right-hand side of \eqref{op.val.} can be for any $\sigma\ge3/2$ estimated as
\begin{eqnarray*}
\mathrm{tr}\left(H_0-\Lambda\right)_-^\sigma \le\,
L_{\sigma,1}^{\mathrm{cl}}\int_{s_0}^\infty\sum_{j=1}^\infty\Big(-\Big(\frac{\pi j}{d(s)} \Big)^2 +W(s)+\Lambda\Big)_+^{\sigma+\frac{1}{2}}\mathrm{d}s +c_2(\Lambda).
\end{eqnarray*}
Combining next the above bound with (\ref{est.}) and using the unitary equivalence between $H_{\Omega_1}$ and $\widetilde{H}_{\widetilde{\Omega}_1}$ we obtain
$$
\mathrm{tr}\left(H_{\Omega_1}-\Lambda\right)_-^\sigma\le\,
L_{\sigma,1}^{\mathrm{cl}}\int_{s_0}^\infty\sum_{j=1}^\infty\Big(-\Big(\frac{\pi j}{d(s)} \Big)^2 +W(s)+\Lambda\Big)_+^{\sigma+\frac{1}{2}}\mathrm{d}s +c_2(\Lambda)
$$
and a simple manipulation of the right-hand side of this inequality leads to
\begin{align}
\mathrm{tr} & \left(H_{\Omega_1} -\Lambda\right)_-^\sigma \nonumber \\[.3em] & \le\,L_{\sigma,1}^{\mathrm{cl}}\int_{\big\{d(s)\ge\pi(W(s)+2\Lambda)^{-1/2}\big\}} \sum_{j=1}^{\left[\frac{1}{\pi}(W(s)+\Lambda)^{1/2} d(s)\right]}\left(W(s)+\Lambda\right)^{\sigma+\frac{1}{2}}\mathrm{d}s+c_2(\Lambda) \nonumber \\[.3em] & \le\frac{L_{\sigma,1}^{\mathrm{cl}}} {\pi}\int_{\big\{d(s)\ge\pi(W(s)+2\Lambda)^{-1/2}\big\}}\left(W(s)+\Lambda\right)^{\sigma+1}\,d(s)\,\mathrm{d}s+c_2(\Lambda) \nonumber \\[.3em] & \label{first result}\le\frac{L_{\sigma,1}^{\mathrm{cl}}}{\pi} \left(\|W\|_{L^\infty(s_0, \infty)}+\Lambda\right)^{\sigma+1} \int_{\big\{d(s)\ge\pi(W(s)+2\Lambda)^{-1/2}\big\}}d(s)\,\mathrm{d}s+c_2(\Lambda).
\end{align}

Now let us pass to the inner part of the spiral associated with operator $H_{\Omega_2}$. We use another version of the `mirroring' trick \cite{F22}.
\begin{figure}[!t]
\begin{center}
\includegraphics[clip, trim=3.5cm 9.5cm 3.5cm 8cm, width=0.75\textwidth]{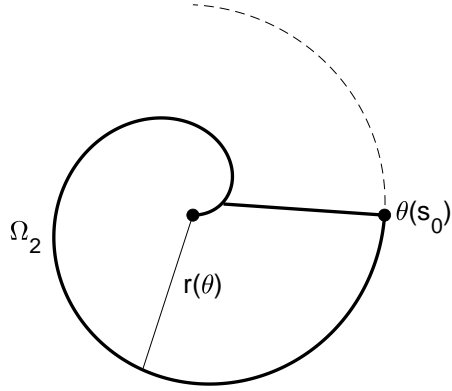}
\end{center}
\caption{The region $\Omega_2$}
\end{figure}
Let $l$ be the segment obtained by extending the straight part of that boundary, that is, the interval $\{s=s_0\}\times\{u\in(0, d(s_0))\}$, to the left up to the boundary of $\Omega_2$. It divides the region into two parts; we denote by $\Omega_2^1$ and $\Omega_2^2$ the upper and lower one, respectively. By Neumann bracketing one gets
\begin{equation}
\label{Omega2}
H_{\Omega_2}\ge H_{\Omega_2^1}\oplus H_{\Omega_2^2},
\end{equation}
where both operators  $H_{\Omega_2^1}$ and $ H_{\Omega_2^2}$ are the restrictions of $H_{\Omega_2}$ to $\Omega_2^1$ and on $\Omega_2^2$, respectively, with the additional Neumann condition imposed on $l$. The spectra of $H_{\Omega_2^{1}}$ and $H_{\Omega_2^{2}}$ are contained in the spectra $H_{\Omega_2^{1,\mathrm{sym}}}$ and $H_{\Omega_2^{2,\mathrm{sym}}}$, respectively, where the latter regions are unions of the former ones and their mirror images with respect to the line spanned by the segment $l$. Using the Berezin inequality \cite{Be72a, Be72b} we estimate the traces of $H_{\Omega_2^{j,\mathrm{sym}}},\: j=1,2,$ as follows
\begin{equation}\label{Berezin}
\mathrm{tr}(H_{\Omega_2^{j,\mathrm{sym}}}-\Lambda)_-^\sigma\le L_{\sigma, 2}^{\mathrm{cl}}\,\mathrm{vol}(\Omega_2^{j,\mathrm{sym}})\,\Lambda^{\sigma+1}, \quad j=1,2,
\end{equation}
where $\Lambda>0$, $\sigma\ge1$ are any numbers and the semiclassical constant
\begin{equation}\label{L2}
L_{\sigma, 2}^{\mathrm{cl}}:= \frac{\Gamma(\sigma+1)}{4\pi\Gamma(\sigma+2)} = \frac{1}{4\pi(\sigma+1)}.
\end{equation}
Finally, combining inequalities (\ref{bracketing}) and (\ref{first result})--(\ref{Berezin}) we arrive at
\begin{align*}
& \mathrm{tr}\left(H_\Omega-\Lambda\right)_-^\sigma \\[.3em] & \le\,\frac{L_{\sigma,1}
^{\mathrm{cl}}}{\pi} \big(\|W\|_{L^\infty(s_0, \infty)}+\Lambda\big)^{\sigma+1} \int_{\big\{d(s)\ge\pi(W(s)+2\Lambda)^{-1/2}\big\}}d(s)\,\mathrm{d}s + c_1 \Lambda^{\sigma+1}+c_2(\Lambda),
\end{align*}
where
\begin{align}
c_1 & := L_{\sigma, 2}^{\mathrm{cl}}\,\big(\mathrm{vol}(\Omega_2^{1, \mathrm{sym}})+\mathrm{vol}(\Omega_2^{2, \mathrm{sym}})\big) \nonumber \\[.3em]
& = 2L_{\sigma, 2}^{\mathrm{cl}}\,\big(\mathrm{vol}(\Omega_2^1)+\mathrm{vol}(\Omega_2^2)\big) =2L_{\sigma, 2}^{\mathrm{cl}}\,\mathrm{vol}(\Omega_2). \label{c}
\end{align}
It remains to check the asymptotic behavior of the quantity \eqref{cc}. Since $\lambda_1(\Lambda)\le\|W\|_{L^\infty(s_0, \infty)}+\Lambda$, it is straightforward to check that
$$
c_2(\Lambda) = \mathcal{O}\Big(\Lambda^2 \int_{\left\{d(s)\ge\frac{\pi}{\sqrt{\Lambda}}\right\}}d(s)\,ds\Big).
$$
holds for for large values of $\Lambda$; this concludes the proof.
\end{proof}

\begin{remark} \label{remark}
{\rm Let us note that Theorem(\ref{Main}) remains valid for smaller powers, $\sigma\ge 1/2$, provided we replace the semiclassical constant $L_{\sigma,1}^{\mathrm{cl}}$ in the right-hand side of \eqref{theorem} by $2r(\sigma,1) L_{\sigma,1}^{\mathrm{cl}}$. On the other hand, we may then set $c_2=0$. This is related with the modification of the Lieb-Thirring inequality for operator-valued potentials defined on the line $\mathbb{R}$ to the powers
$\sigma\ge 1/2$ and the constant $r(\sigma,1) L_{\sigma,1}^{\mathrm{cl}}$ where $r(\sigma,1)\le 2$ if $\sigma<3/2$ \cite{ELW04} which we used to derive the estimate \eqref{count.}. Indeed, for $\sigma\ge1/2$ we have
\begin{align}
\mathrm{tr} & \,\Big(-\frac{\partial^2}{\partial\,s^2}\otimes \,I_{L^2(\mathbb{R})}+ L(s,\hat{W},\Lambda')\Big)_-^\sigma \nonumber \\[.3em] &
\le\,r(\sigma,1)\,L_{\sigma,1}^{\mathrm{cl}}\int_{\mathbb{R}}\sum_{j=1}^\infty\Big(-\Big(\frac{\pi j}{\hat{d}(s)} \Big)^2 +\hat{W}(s)+\Lambda'\Big)_+^{\sigma+\frac{1}{2}}\mathrm{d}s
\nonumber \\[.3em] & =\,2 r(\sigma,1)\,L_{\sigma,1}^{\mathrm{cl}}\int_{s_1}^\infty\sum_{j=1}^\infty\Big(-\Big(\frac{\pi j}{d(s)} \Big)^2 +W(s)+\Lambda'\Big)_+^{\sigma+\frac{1}{2}}\mathrm{d}s \nonumber \\[.3em] & \label{ELW}\le\frac{2r(\sigma,1)\,L_{\sigma,1} ^{\mathrm{cl}}}{\pi} \big(\|W\|_{L^\infty(s_0, \infty)}+\Lambda'\big)^{\sigma+1} \int_{\big\{d(s)\ge\pi(W(s)+2\Lambda)^{-1/2}\big\}}d(s)\,\mathrm{d}s.
\end{align}
This completes the argument leading to \eqref{count.} and at the same time, allows us to get the result mentioned in the opening of the remark by repeating the steps of the previous proof.
}
\end{remark}

\section{Concluding remarks} \label{s: conl}
\setcounter{equation}{0}

\subsection{Optimality of the bound (\ref{theorem})}

Let us now show that the bound of Theorem~\ref{Main} is asymptotically sharp in the sense that the dependence on $\Lambda$ the leading term in \eqref{asy_corollary} cannot be improved. To provide an example proving this claim, consider a spiral region $\Omega$ such that its parallel coordinates representation (\ref{parametrization}) outside a compact area satisfies $d(s)=s^{-1}$ starting from some $s_0>0$. As indicated in the discussion of assumption~\eqref{assa}, one can choose $s_0$ in such a way that $w:=\|d(s)\gamma(s)\|_{L^\infty(s_0,\infty)}<1$.

Using Dirichlet bracketing method one estimates $H_\Omega$ from above as follows,
\begin{equation}\label{Dir.brack.}
H_\Omega\le H_1\oplus H_2,
\end{equation}
where $H_1$ and $H_2$ are the Dirichlet restrictions of $H_\Omega$ to $\Omega_1\subset\Omega$ satisfying $s>s_0$ and on $\Omega_2=\Omega\setminus\bar{\Omega}_1$. Since $\Omega_2$ is bounded, its contribution to the eigenvalue moment count is of the standard form \cite{W12},
\begin{equation} \label{Omega2}
\mathrm{tr}\left(H_2-\Lambda\right)_-^\sigma = L^{\mathrm{cl}}_{\sigma,2}\,
\mathrm{vol}(\Omega_2)\,\Lambda^{\sigma+1}+ \bar{o}(\Lambda^{\sigma+1}), \quad\sigma\ge0,\; \Lambda\to\infty,
\end{equation}
with $L_{\sigma, 2}^{\mathrm{cl}}$ given by \eqref{L2}. To deal with $H_1$ we pass to the unitary equivalent operator \eqref{unitary eq.} acting as
$$
\widetilde{H}=-\frac{\partial}{\partial s}\left(\frac{1}{(1-u\gamma(s))^2}\frac{\partial}{\partial s}\right)-\frac{\partial^2}{\partial u^2}+\widetilde{W}(s, u)
$$
on $\widetilde{\Omega}_1:=\{(s, u): s>s_0, u\in(0, d(s))\}$ with the effective potential $\widetilde{W}$ given by (\ref{tildeW}). It is straightforward to check that
$$
\widetilde{H}\le \frac{1}{(1-w)^2}(-\Delta_\mathrm{D})+\|\widetilde{W}\|_{L^\infty(\widetilde{\Omega}_1)},
$$
where $\Delta_\mathrm{D}$ is the Dirichlet Laplacian in $L^2(\widetilde{\Omega}_1)$. Consequently, for any $\sigma\ge0$ we have
\begin{equation}\label{Berg}
\mathrm{tr}(\widetilde{H}-\Lambda)_-^\sigma\ge\frac{1}{(1-w)^{2\sigma}}
\mathrm{tr}\left(-\Delta_\mathrm{D}-(1-w)^2(\Lambda-\|\widetilde{W}\|_{L^\infty(s>s_0, u\in(0, d(s))})\right)_-^\sigma.
\end{equation}
The right-hand side can be estimated using the asymptotic properties of the spectral counting function of the Dirichlet Laplacian on horn-shaped regions \cite{V92}:
\begin{theorem}
Let functions $f_j: [0, \infty)\to \mathbb{R}_+,\, j=1,2,$ be right-continuous and decreasing to zero. Consider the region
$$
D=\{(s, u): s>0, -f_1(s)<u< f_2(s)\} \subset \mathbb{R}^2
$$
and suppose that $f(s):=f_1(s)+f_2(s),\, s>0,$ satisfies
$$
\int_0^\infty e^{-t f(s)^{-2}}\,\mathrm{d}s< \infty,\quad t > 0.
$$
Then for the number $N_D(\lambda)$ of eigenvalues of the Dirichlet Laplacian on $D$ that are less than $\lambda$ we have
\begin{equation}\label{count.}
N_D(\lambda)\sim\int_0^\infty \sum_{k=1}^\infty\Big(\Big(\frac{\lambda}{\pi^2}-\frac{k^2}{f^2(s)}\Big)_+\Big)^{1/2}\,\mathrm{d}s,\quad \lambda\to\infty,
\end{equation}
where $f(t)\sim g(t)$ means that $f(t)/g(t)\to 1$ as $t\to\infty$.
\end{theorem}

In our case $f_1=0, \,f_2(s)=f(s)=d(s)$ and $\lambda=\Lambda$; a series of simple estimates then gives
\begin{align*}
\int_0^\infty  & \sum_{k=1}^\infty\Big(\Big(\frac{\Lambda}{\pi^2}-\frac{k^2}{d^2(s)}\Big)_+\Big)^{1/2}\,\mathrm{d}s \\[.3em]
& =\int_0^\infty \sum_{\left\{k\le \frac{\sqrt{\Lambda}d(s)}{\pi}\right\}}\Big(\Big(\frac{\Lambda}{\pi^2}-\frac{k^2}{d^2(s)}\Big)_+\Big)^{1/2}\,\mathrm{d}s \\[.3em] & =\int_{\left\{d(s)\ge\frac{\pi}{\sqrt{\Lambda}}\right\}}^\infty \sum_{\left\{k\le \frac{\sqrt{\Lambda} d(s)}{\pi}\right\}}\Big(\Big(\frac{\Lambda}{\pi^2}-\frac{k^2}{d^2(s)}\Big)\Big)^{1/2}\,\mathrm{d}s \\[.3em] & \ge\int_{\left\{d(s)\ge\frac{\pi}{\sqrt{\Lambda}}\right\}}^\infty \sum_{\left\{k\le \frac{\sqrt{\Lambda} d(s)}{\sqrt{2}\pi}\right\}}\Big(\Big(\frac{\Lambda}{\pi^2}-\frac{k^2}{d^2(s)}\Big)\Big)^{1/2}\,\mathrm{d}s \\[.3em] & \ge\frac{\Lambda}{2\pi^2}\int_{\left\{d(s)\ge\frac{\pi}{\sqrt{\Lambda}}\right\}}d(s)\,\mathrm{d}s,
\end{align*}
where in the last step we simply took the lower bound to the square root expression times the number of summands. Using next the fact that we have $d(s)=s^{-1}$ in combination with \eqref{Berg} in which we divide the potential into the sum of two equal parts, we get
\begin{align}
\nonumber \mathrm{tr}&(\widetilde{H}-\Lambda)_-^\sigma \ge\frac{1}{(1-w)^{2\sigma}}\left(\frac12 (1-w)^2\big( \Lambda-\|\widetilde{W}\|_{L^\infty(s>s_0, u\in(0, d(s))}\big)\right)^\sigma \\[.3em] & \nonumber \quad\times N_{\widetilde{\Omega}_1}\left(\frac12 (1-w)^2 \big(\Lambda-\|\widetilde{W}\|_{L^\infty(s>s_0, u\in(0, d(s))}\big)\right)
\\[.3em] & \nonumber >\left(\frac12 \big(\Lambda-\|\widetilde{W}\|_{L^\infty(s>s_0, u\in(0, d(s))}\big)\right)^\sigma \frac{1}{4\pi^2}(1-w)^2 \big(\Lambda-\|\widetilde{W}\|_{L^\infty(s>s_0, u\in(0, d(s))}\big)
\\[.3em] & \nonumber \quad\times\int_{\left\{d(s)\ge\pi\sqrt{2} (1-w)^{-1} \big(\Lambda-\|\widetilde{W}\|_{L^\infty(s>s_0, u\in(0, d(s)))}\big)^{-1/2}\right\}} d(s)\,\mathrm{d}s \\[.3em] & \label{end}
=\frac{(1-w)^2}{2^{\sigma+3}\pi^2}\,\Lambda^{\sigma+1} \ln\Lambda\,(1+\bar{o}(1)),
\end{align}
which differs from \eqref{asy_corollary} by a multiplicative constant only. We note that $w$ in the above discussion can be chosen arbitrarily small and the ratio of the two constants,
\begin{equation} \label{ratio}
2^{\sigma+3}\pi\,L_{\sigma,1}^{\mathrm{cl}} = 2^{\sigma+2}\, \sqrt{\pi}\, \frac{\Gamma(\sigma+1)}{\sqrt{4\pi} \Gamma(\sigma+\frac32)},
\end{equation}
should not be smaller than one. In fact, the bound \eqref{end} illustrates that the estimate in \eqref{asy_corollary} is getting worse with increasing $\sigma$, since we have \cite{TE52}
$$
\frac{\Gamma(z)}{\Gamma(z+1/2)}=\frac{1}{\sqrt{z}}(1+\bar{o}(1))\quad\text{as}\quad z\to\infty.
$$
For half-integer values, $\sigma=n+\frac12$ with $n\in\mathbb{N}_0$ the ratio \eqref{ratio} can be expressed explicitly as
$$
2^{\sigma+3}\pi\,L_{\sigma,1}^{\mathrm{cl}}= 2^{3/2}\,\sqrt{\pi}\, \frac{(2n+1)!!}{(n+1)!},
$$
in particular, for $\sigma=3/2$ the value is $3\sqrt{2}\pi \approx 13.3$.

\subsection{Multi-arm spirals}

Let us note finally that the obtained result extends easily to the situation where the region is determined by a multi-armed spiral.

Let $\Gamma_0=(r_0(\theta), \theta)$ be a shrinking spiral satisfying the assumptions of Theorem(\ref{Main}) and let $\Gamma_m$ be the union of $m$ `angularly shifted' spirals $\Gamma_j=(r_j(\theta), \theta)$, where
$r_j(\theta):=r\left(\theta-\theta_j\right)$ for $0\le j\le m-1$ corresponding to the partition $0=\theta_0<\ldots<\theta_{m-1}<2\pi$. Let $H_{\Omega_{\Gamma_m}}$ be the Laplace operator defined on $\Omega_{\Gamma_m}=\mathbb{R}^2\setminus\Gamma_m$ with the Dirichlet conditions imposed on $\Gamma_m$. It is obvious that $H_{\Omega_{\Gamma_m}}$ decomposes into the direct sum of $m$ Dirichlet Laplacians unitarily equivalent to operators $H_j$ acting on the `rotated' spiral-shaped domains $\Omega_j:=\{(r, \theta): r\in (\mathrm{max}\{0, r(\theta-\theta_{j-1})\}, r(\theta-\theta_j)\}$, and consequently, it is enough to estimate separately $\mathrm{tr}\left(H_j-\Lambda\right)_-^\sigma,\, 0\le j\le m-1$. By a straightforward modification of Theorem~(\ref{Main}) we get
\begin{align}
\nonumber\mathrm{tr}\left(H_j-\Lambda\right)_-^\sigma\le & \frac{L_{\sigma,1} ^{\mathrm{cl}}}{\pi} \left(\|W_j\|_{L^\infty(s'_{1, j}, \infty)}+\Lambda\right)^{\sigma+1} \int_{\left\{d_j(s)\ge \pi (W_j(s)+\Lambda)^{-1/2}\right\}}d_j(s)\,\mathrm{d}s \\[.3em] \label{multi-arm} & + c_j\Lambda^{\sigma+1} +c_2^j(\Lambda),
\end{align}
where $d_j(s), W_j(s), s'_{0, j}, c_j$, and $c_2^j(\Lambda)$ are natural modifications of the quantities appearing in the main result. To estimate $\mathrm{tr}\left(H_{\Omega_{\Gamma_m}}-\Lambda\right)_-^\sigma$ one has to sum the expressions on the right-hand side of \eqref{multi-arm}. In particular, in the asymptotic regime, $\Lambda\to\infty$, we have
$$
\mathrm{tr}\left(H_{\Omega_{\Gamma_m}}-\Lambda\right)_-^\sigma \le  \Lambda^{\sigma+1}\bigg(\frac{L_{\sigma,1}
^{\mathrm{cl}}}{\pi}\sum_{j=0}^{m-1}\int_{\left\{d_j(s)\ge\frac{\pi}{\sqrt{\Lambda}}\right\}}d_j(s)\,\mathrm{d}s+\tilde{c} m\bigg)(1+\bar{o}(1)),
$$
where the constant $\tilde{c}:=\underset{0 \leq j \leq m-1}{\max}\,c_j$.

\section*{Appendix}

Here we provide proof of Theorem~\ref{operator valued1} which was skipped in Sec.~\ref{Main result}. Let $\lambda_1, \ldots, \lambda_N$ be the negative eigenvalues of $-\frac{\partial^2}{\partial\,s^2}\otimes \,I_{L^2(\mathbb{R})}+ Q(s)$ and denote by $P_N$ the projection on the linear span of the corresponding eigenfunctions. It is easy to see that
\begin{equation}\label{span}
\mathrm{tr}\left(-\frac{\partial^2}{\partial\,s^2}\otimes \,I_{L^2(\mathbb{R})}+ Q(s)\right)_-^\sigma\le \mathrm{tr}(H_N)_-^\sigma,\quad\sigma\ge0,
\end{equation}
holds for the finite-dimensional restriction of the operator,
$$H_N=P_N\,\left(-\frac{\partial^2}{\partial\,s^2}\otimes \,I_{L^2(\mathbb{R})}+ Q(s)\right)\,P_N.
$$
The expression on the right-hand side of \eqref{span} is nothing but the Riesz mean of the order $\sigma$ of the negative eigenvalues of the $N\times N$-system of ordinary differential equations $-\frac{\partial^2}{\partial\,s^2}\otimes\,I+P_N Q(s) P_N$ acting on $L^2(\mathbb{R}_+, \mathbb{C}^N)$ with the Neumann boundary condition at $s=0$, where $I$ denotes the identity operator on $\mathbb{C}^N$. This allows us to use the Lieb-Thirring inequality for general second-order differential operators with matrix-valued potentials on the positive half-line with Neumann condition at the origin proved in \cite{M15}. In this way we get
\begin{align*}
\mathrm{tr} (H_N)_-^\sigma & \le L_{\sigma, 1}^{\mathrm{cl}}\int_0^\infty\mathrm{tr}\,\left(P_N Q(s) P_N\right)^{\sigma+1/2}(s)\,\mathrm{d}s+ \frac12\,N \lambda_1^{(N)} \\[.3em] & \le
 L_{\sigma, 1}^{\mathrm{cl}}\int_0^\infty\mathrm{tr}\,Q^{\sigma+1/2}(s)\,\mathrm{d}s+\frac12\,N \lambda_1^{(N)}
 \end{align*}
with $\lambda_1^{(N)}$ given by \ref{Nlambda} which is the sought result.

\subsection*{Acknowledgements}

The work D.B. has been supported by the Czech Science Foundation (GA\v{C}R) within the project 21-07129S and by the Czech-Polish project BPI/PST/2021/1/00031.
The work P.E. has been supported by the Czech Science Foundation (GA\v{C}R) within the project 21-07129S and the EU project CZ.02.1.01/0.0/0.0/16\textunderscore 019/0000778. The authors are obliged to Rupert Frank for a useful discussion.


\end{document}